\documentclass[12pt]{article}
\usepackage{amssymb,amsmath}
\usepackage{graphicx}
\usepackage{amsfonts}
\usepackage{float}
\usepackage{flafter}
\usepackage{color}

\labelwidth\leftmargini\advance\labelwidth-\labelsep

\def\R{\relax\ifmmode I\!\!R\else$I\!\!R$\fi}

\def\Z{\relax\ifmmode Z\!\!\!Z\else$Z\!\!\!Z$\fi}

\def\C{\relax\ifmmode C\!\!\!\!I\else$C\!\!\!\!I$\fi}

\def\K{\relax\ifmmode I\!\!K\else$I\!\!K$\fi}

\def\N{\relax\ifmmode I\!\!N\else$I\!\!N$\fi}

\newcounter{defcounter}[section]
{\vspace{0.1cm}\begin{sloppypar}\noindent\stepcounter{defcounter}{\bfseries
Definition
      \thesection.\thedefcounter}}%
{\end{sloppypar}\vspace{0.1cm}}

\newtheorem{theorem}{Theorem}[section]

\newtheorem{proposition}{Proposition}[section]

\newcommand{\proof}{{\bf Proof.} }

\newcommand{\qed}{\hfill $\square$}

\begin{document}
\thispagestyle{empty}
\begin{center}
{\Large {\bf Continued logarithm representation of real numbers}}
\end{center}
\begin{center}J\"org Neunh\"auserer\\
Leuphana University L\"uneburg, Germany\\
joerg.neunhaeuserer@web.de
\end{center}
\begin{center}
\begin{abstract} We introduce the continued logarithm representation of real numbers and prove results on the occurrence and frequency of digits with respect to this representation.
. ~
\end{abstract}
\end{center}
\section{Introduction}
The continued fraction representation of real numbers is intensively studied in number theory, see \cite{[CO]} and references there in. In this paper we consider the continued logarithm representation of real numbers, which we now introduce. For $m\ge 3$ and a sequence $(d_{k})\in\{1,\dots,m-1\}^{\mathbb{N}}$ the continued logarithm to base $m$ is given by
\[ [(d_{k})]_{m}=\lim_{k\to\infty}\log_{m}(d_{n}+\log_{m}(d_{n-1}+\log_{m}(\dots+\log_{m}(d_{1})\dots)\]
The limit exists since the maps
\[ T_{d}(x)=\log_{m}(d+x)\]
are contractions on $[0,1]$ for $d\in \{1,\dots,m-1\}$. For a finite sequence $(d_{k})\in\{1,\dots,m-1\}^{n}$ we
consider the iterated functions
\[ [(d_{k})]_{m}(x)=T_{d_{n}}\circ\dots\circ T_{d_{1}}(x),\]
which map $[0,1]$ to nested subintervals. In this case $[(d_{k})]_{m}$ denotes the closed interval
 \[ [(d_{k})]_{m}=[~[(d_{k})]_{m}(0)~,~[(d_{k})]_{m}(1)~].\]
Obviously for an infinite sequence $(d_{k})$ we have
$[(d_{k})]_{m}\in [(d_{1},\dots,d_{n})]_{m}$
for all $n\in\mathbb{N}$.  As expected we have obtain the following theorem on the continued logarithm representation.
\begin{theorem}
All real numbers in the interval $[0,1]$ have a continued logarithm representation to base $m$ and except to a countable set this representation is unique.
\end{theorem}
\proof Consider the map $f(x)=m^x$ mod $1$. For $x\in [0,1)$ let $d_{k}=i$ if $f^{k-1}(x)\in [\log_{m}(i),\log_{m}(i+1))$ for $i\in\{1,\dots,m-1\}$.
Since the maps $T_{d}(x)$ are the inverse branches of $f$ the point $x$ is contained in the interval $[d_{1},\dots,d_{n}]$ for all sequences $(d_{1},\dots,d_{n})$ and hence $[(d_{k})]_{m}=x$. Thus we have constructed a continued logarithm representation for all $x\in[0,1)$. The representation of $1$ is obviously $(m-1)$. Furthermore the interior of the intervals $[(d_{1},\dots ,d_{n}]_{m}$ are disjoint for all $(d_{k})\in\{1,\dots,m-1\}^{n}$. Hence the representation of $x\in[0,1]$ is ambiguous if and only if $x=[(d_{k})]_{m}(1)$ for some sequence $(d_{k})\in\{1,\dots,m-1\}^{n}$. But the set of these sequences is countable.
\qed~\\~\\
As far as we know no results on the continued logarithm representation were published. In the following we will find results on the occurrence of digits, sets with restricted digits and the frequency of digits with respect to this this representation.
\section{Occurrence of digits}
As in the case of the usual power-series representation of real numbers to base $b\ge 2$ the following result on the occurrence of digits holds in the continued logarithm representation.
\begin{theorem}
In the continued logarithm representation to base $m\ge 3$ of almost all real numbers in $[0,1]$ all digits $i\in \{1,\dots,m-1\}$ occur infinitely many times.
\end{theorem}
\proof In the following $|I|$ denotes the Lebesgue measure of $I\subset [0,1]$. Fix $m\ge 3$. Since the maps $T_{d}(x):[0,1]\to[0,1]$ are conformal contractions with $T_{i}(0,1)\cap T_{j}(0,1)=\emptyset$ for $i\not=j$, there are contraction constants $c_{d}\in(0,1)$ with $\sum_{d\in D }c_{d}=1$ such that
\[ |[d_{1},\dots, d_{n}]_{m}|\le C\prod_{i=1}^{n} c_{d_{i}} \]
where $C>0$ is a constant independent of $n\in\mathbb{N}$. Let $R_{l}(\tilde D)$ be the set of all real numbers in $[0,1]$ where the first $l$ digits in the continued logarithm representation come from $D$ and the other digits come a proper subset $\tilde D\subset D$. We have
\[R_{l}(\tilde D)\subseteq\bigcup_{d_{1},\dots d_{l}\in D, d_{l+1},\dots,d_{n}\in\tilde D }[(d_{1},\dots d_{l},d_{l+1},\dots,d_{n})]_{m}\]
for all $n>l$. Hence we obtain
\[ |R_{l}(\tilde D)|\le \sum_{d_{1},\dots, d_{l}\in D, d_{l+1},\dots,d_{n}\in\tilde D }|~[(d_{1},\dots ,d_{l},d_{l+1},\dots,d_{n})]_{m}~|\]
\[ \le \sum_{d_{1},\dots ,d_{l}\in D, d_{l+1},\dots,d_{n}\in\tilde D }C\cdot c_{d_{1}}\cdot\dots \cdot c_{d_{l}}\cdot c_{d_{l+1}}\cdot\dots\cdot c_{d_{n}}\]
\[ = C\cdot(\sum_{d\in D} c_{d})^{l}\cdot(\sum_{d\in \tilde D}c_{d})^{n-l}\]
for all $n>l$. Since $\sum_{d\in \tilde D}c_{d}<1$ we have $|R_{l}(\tilde D)|=0$.
Now consider the set of all real numbers in $[0,1]$ for which not all digits $i\in \{1,\dots,m-1\}$ occur infinitely many times. This is the set
\[ \bigcup_{i\in\{1,\dots,m-1\}}\bigcup_{n\in\mathbb{N}_{0}}\bigcup_{d_{1},\dots,d_{n}\in\{1,\dots,m-1\}}\{[(d_{k})]_{m}~|~d_{k}\not=i~\forall k>n\}.\]
The Lebesgue measure of this set is zero since the countable union of sets with Lebesgue measure zero has Lebesgue measure zero as well. Taking the complement gives the result.
 \qed~\\~\\
Let $m\ge 4$ and let $D\subset\{1,\dots,m-1\}$ be a subset with more than one Element. We consider the set $[D^{\mathbb{N}}]_{m}$ of all reel numbers in $[0,1]$ that have a continued logarithm representation to base $m$ with digits in $D$.  This set is obviously uncountable and from the proof of the last theorem we know that it is totally disconnected. It is natural to ask for the Hausdorff dimension $\dim_{H}[D^{\mathbb{N}}]_{m}$ of this set. We refer to \cite{[FA]} or \cite{[PE]} for an introduction to dimension theory. We will estimate the Hausdorff dimension of $[D^{\mathbb{N}}]_{m}$ using the following well know theorem:
\begin{theorem}
Let $T_{i}:\mathbb{R}^{n}\to \mathbb{R}^{n}$ for $i=1,\dots, r$ be a iterated function system consisting of conformal contractions fulfilling the open set condition, which means that there is an open set $O\subseteq\mathbb{R}^{n}$ such that $T_{i}(O)\subseteq O$ and $T_{i}(O)\cap T_{j}(O)=\emptyset$ for $i\not=j$. If we have
\[ c_{i}<|T'_{i}(x)|<C_{i}\]
on $O$ and $d,D>0$ are given by
\[\sum_{i=1}^{r}c_{i}^{d}=1\qquad\sum_{i=1}^{r}C_{i}^{D}=1\]
the Hausdorff dimension of the unique compact set $K$ with $K=T_{1}(K)\cup\dots T_{r}(K)$ is bounded by
\[ d<\dim_{H}K<D,\]
\end{theorem}
By theorem 9.9 of \cite{[FA]} we immediately get the upper bound in this theorem and the lower bound follows from theorem 3.15 of \cite{[MU]}, which is in fact more general. We now obtain:
\begin{theorem} For $D\subset \{1,\dots,m-1\}$
we have
\[L_{n}\le \dim_{H} [D^{\mathbb{N}}]_{m}\le U_{n} \]
for all $n\ge1$, where $L_{n}$ and $U_{n}$ are given by
\[ \sum_{d_{1},\dots ,d_{n}\in D} [(d_{k})]'(0)^{L_{n}}=1\quad \sum_{d_{1},\dots ,d_{n}\in D} [(d_{k})]'(1)^{U_{n}}=1\]
\end{theorem}
\proof The set $[D^{\mathbb{N}}]_{m}$ is the attractor of the iterated function system
\[ \{[(d_{k})]_{m}(x)~|~d_{1},\dots ,d_{n}\in D\}\]
for all $n\ge 1$, see \cite{[HU]} or chapter nine of \cite{[FA]}. Note that the iterated function system is conformal and fulfils the open set condition since the images of the maps may intersect only in the boundary,  see \cite{[MU]}. Furthermore note that
\[ \max\{[(d_{k})]'_{m}(x)~|~x\in[0,1]\}=[(d_{k})]'_{m}(0) \]
\[ \min\{[(d_{k})]'_{m}(x)~|~x\in[0,1]\}=[(d_{k})]'_{m}(1), \]
which means the the contraction rates of the maps that generate the iterated function system are bounded from above and below.
Our result thus follows from theorem 2.2.
\qed~\\~\\
We use Mathematica to solve the equations in the last proposition for $m=4$ and obtain
\[ \dim_{H} [\{1,2\}^{\mathbb{N}}]_{4}=0.81\pm 0.01\]
\[ \dim_{H} [\{1,3\}^{\mathbb{N}}]_{4}=0.66\pm 0.01\]
\[ \dim_{H} [\{2,3\}^{\mathbb{N}}]_{4}=0.45\pm 0.01\]
Compare this with the classical result of Hausdorff \cite{[HA]} that the dimension of the set of real numbers with one delted digit in the powers series representation to base $3$ is $\log(2)/\log(3)=0.630\dots$, no matter which digit is deleted.
\section{Frequency of digits}
Let $\mathfrak{f}_{i}([(d_{k})]_{m})$ be the frequency of the digit $i\in\{1,\dots,m-1\}$ in the continued logarithm representation $[(d_{k})]_{m}$ of a real number in $[0,1]$, that is
\[ \mathfrak{f}_{i}([(d_{k})]_{m})=\lim_{n\to\infty}\frac{\sharp\{k\in\{1,\dots,n\}|d_{k}=i\}}{n},\]
provided that the limit exist. We consider sets with of real numbers with given frequencies of the continued logarithm representation to base $m$. For a probability vector $(p_{1},\dots,p_{m-1})\in(0,1)^{m-1}$ let
\[ \mathfrak{F}_{m}(p_{1},\dots,p_{m-1})=\{[(d_{k})]_{m}\in[0,1]~|~\mathfrak{f}_{i}([(d_{k})]_{m})=p_{i},~i=1,\dots,m-1\}\]
We first prove an upper bound on the Hausdorff dimenison of theses sets
\begin{proposition}
\[ \dim_{H}\mathfrak{F}_{m}(p_{1},\dots,p_{m-1})\le \frac{-\sum_{i=1}^{m-1}p_{i}\log(p_{i})}{\sum_{i=1}^{m-1}p_{i}\log(\log(m-1)+\log(m)i)}\]
\end{proposition}
\proof
We will prove the dimension estimate for $T_{m-1}(\mathfrak{F}_{m}(p_{1},\dots,p_{m-1}) )$. The result follows since $T^{-1}_{m-1}x=m^x-m+1$ is Lipschitz on $[\log_{m}(m-1),1]$ and hence does not increase Hausdorff dimension, see corollary 2.4 of \cite{[FA]}. \\
Again we write $[(d_{1},\dots,d_{n})]_{m}$ for the interval $[[(d_{1},\dots,d_{n})]_{m}(0),[(d_{1},\dots,d_{n})]_{m}(1)]$ and denote the length of an interval  $I\subseteq\mathbb{R}$ by $|I|$. For
$[(d_{k})]_{m}\in T_{m-1}(\mathfrak{F}_{m}(p_{1},\dots,p_{m-1}))$ we have
\[ |[(d_{1},\dots,d_{n})]_{m}|\le \max\{(T_{d_{1}}\circ\dots\circ T_{d_{n}})'(x)|x\in [\log_{m}(m-1),1]\}\]
\[ \le \prod_{i=1}^{n}\frac{1}{\log(m)}\max\{\frac{1}{x+d_{i}}|x\in [\log_{m}(m-1),1)\}=\prod_{i=1}^{n}\frac{1}{\log(m)(\log_{m}(m-1)+d_{i})}\]\[
=(\prod_{i=1}^{n}\log(m-1)+\log(m)d_{i}))^{-1},
\]
hence
\[\liminf_{n\to\infty}-\frac{1}{n}\log(|[(d_{1},\dots,d_{n})]_{m}|)\]
\[ \ge \liminf_{n\to\infty}\frac{1}{n}\sum_{i=1}^{n}\log(\log(m-1)+\log(m)d_{i}))\]
\[ =\sum_{i=1}^{m-1}p_{i}\log(\log(m-1)+\log(m)i).\]
In the last equation we use the frequency of digits in $[(d_{k})]_{m}$. Now consider a Borel probability measure on $[0,1]$ with \[ \mu([(d_{1},\dots,d_{n})]_{m})=\prod_{i=1}^{n}p_{d_{i}}.\]
We obviously have
\[ \lim_{n\to\infty}\frac{1}{n}\log\mu([(d_{1},\dots,d_{n})]_{m})=-\sum_{i=1}^{m-1}p_{i}\log(p_{i})\]
and hence
\[ \limsup_{n\to\infty}\frac{\log\mu([(d_{1},\dots,d_{n})]_{m})}{\log(|[(d_{1},\dots,d_{n})]_{m}|)}\]\[
\le \frac{-\sum_{i=1}^{m-1}p_{i}\log(p_{i})}{\sum_{i=1}^{m-1}p_{i}\log(\log(m-1)+\log(m)i)}=:U_{m}(p_{1},\dots,p_{m-1})\]
for all $[(d_{k})]_{m}\in T_{m-1}(\mathfrak{F}_{m}(p_{1},\dots,p_{m-1}))$. Note that the intervals $[(d_{1},\dots,d_{n})]_{m}|$ constitute a nested
sequence of partitions with
\[ c_{1}^{n}<|[(d_{1},\dots,d_{n})]_{m}|<c_{2}^{n}.\]
Thus we obtain
\[ \liminf_{\epsilon\to 0}\frac{\log\mu((x-\epsilon,x+\epsilon))}{\log(\epsilon)}\le U_{m}(p_{1},\dots,p_{m-1}) \]
for all $x\in T_{m-1}(\mathfrak{F}_{m}(p_{1},\dots,p_{m-1}))$.
This means that the lower local dimension of the measure $\mu$ is bounded by $U_{m}$ and by theorem 7.2 of \cite{[PE]} we obtain
\[ \dim_{H}T_{m-1}(\mathfrak{N}_{m})\le U_{m}(p_{1},\dots,p_{m-1}) \]
\qed ~\\~\\
From the last proposition we obtain a theorem which is striking compared with Borel's \cite{[BO]} classical result that almost all real numbers are normal with respect to usual powers series representations.
\begin{theorem}
For all $m\ge 3$ the set of real numbers in $[0,1]$ that have continued logarithm representation to base $m$ with given  frequencies has Hausdorff dimension less than one.
\end{theorem}
\proof
Let $d>0$ be the solution of
\[ \sum_{i=1}^{m-1}(\log(m-1)+\log(m)i)^{-d}=1 \]
and let $(p_{i})=((\log(m-1)+\log(m)i)^{-d})$ be the corresponding probability vector. The function $U_{m}(p_{1},\dots,p_{m-1})$ attains its maximum for this probability vector and the value of the maximum is $d$. Now observe that
\[\sum_{i=1}^{m-1}(\log(m-1)+\log(m)i)^{-1}<1\]
for $m\ge 3$ hence $d<1$, which completes the proof.
\qed~\\~\\
In the case $m=3$ we have
\[ \dim_{H}\mathfrak{F}_{3}(p,1-p)\le \frac{-p\log(p)-(1-p)\log(1-p)}{p\log(\log(2)+\log(3))+(1-p)\log(\log(2)+2\log(3))}\]
The graph of the upper bound is displayed below

We conjecture that $\dim_{H}\mathfrak{F}_{3}(p,1-p)$ is in fact an unimodal function, but to find an explicit expression for this function seems to be quite difficult.


\begin{thebibliography}
\small
\bibitem{[BO]} E. Borel, Les probabilit´es d´enombrables et leurs applications arithm´etiques. Rend.
Circ. Mat. Palermo 27, 247-271, 1909.
\bibitem{[FA]} K.J. Falconer, Fractal Geometry - mathematical foundations
and applications, Wiley, New York, 1990.
\bibitem{[HA]} F. Hausdorff, Dimension und \"ausseres Mass, Math. Annalen 79, 157-179, 1919.
\bibitem{[CO]} D. Hensley, Continued Fractions, World Scientific Pub Co., New Jersey, 2006.
\bibitem{[HU]} J.J. Hutchinson, Fractals and self-similarity, Indiana
Univ. Math. J. 30, 271-280, 1981.
\bibitem{[MU]}  R.D.Mauldin, M.Urbanski, Dimensions and measures in iterated function systems,Proc. London Math. Soc.(3) 73, 1996.
\bibitem{[PE]} Ya. Pesin, Dimension Theory in Dynamical Systems - contemporary views and ap-
plications, University of Chicago Press, Chicago, 1997.
\end{thebibliography}
\end{document}